\documentclass[10pt]{amsart}

\usepackage{latexsym, amsmath,amssymb}
\usepackage{color}

\renewcommand{\epsilon}{\varepsilon}
\newcommand{\newsection}[1]
{\subsection{#1}\setcounter{theorem}{0} \setcounter{equation}{0}
\par\noindent}

\newtheorem{theorem}{Theorem}

\newtheorem{lemma}[theorem]{Lemma}
\newtheorem{corr}[theorem]{Corollary}

\newtheorem{proposition}[theorem]{Proposition}
\newtheorem{deff}[theorem]{Definition}

\newcommand{\bth}{\begin{theorem}}
\newcommand{\ble}{\begin{lemma}}
\newcommand{\bcor}{\begin{corr}}

\newcommand{\bdeff}{\begin{deff}}

\newcommand{\bprop}{\begin{proposition}}
\newcommand{\ele}{\end{lemma}}
\newcommand{\ecor}{\end{corr}}
\newcommand{\edeff}{\end{deff}}

\newcommand{\eprop}{\end{proposition}}

\renewcommand{\l}{\lambda}

\newcommand{\supp}{\text{supp}\ }
\renewcommand{\Pi}{\varPi}

\renewcommand{\epsilon}{\varepsilon}

\newcommand{\K}{{\mathcal K}}
\newcommand{\R}{{\mathbb R}}

\newcommand{\la}{{\langle}}
\newcommand{\ra}{{\rangle}}
\newcommand{\cd}{{\,\cdot\,}}
\newcommand{\bdy}{{\partial\K}}
\newcommand{\ext}{{\R^3\backslash\K}}
\renewcommand{\S}{{\mathbb{S}}}

\begin{document}

\title[Long-time existence of quasilinear wave equations]
{
The lifespan for 3-dimensional quasilinear wave equations in exterior domains
}

\author{John Helms}

\author{Jason Metcalfe}
\address{Department of Mathematics, University of North Carolina,
  Chapel Hill, NC  27599-3250}
\email{johelms@email.unc.edu, metcalfe@email.unc.edu}

\thanks{The authors were supported in part by NSF grant DMS-1054289.
The second author was additionally supported by NSF grant DMS-0800678}

\begin{abstract}
This article focuses on long-time existence for quasilinear wave
equations with small initial data in exterior domains.  The
nonlinearity is permitted to fully depend on the solution at the
quadratic level, rather than just the first and second derivatives of
the solution.  The corresponding lifespan bound in the boundaryless
case is due to Lindblad, and Du and Zhou first proved such long-time
existence exterior to star-shaped obstacles.  Here we relax the
hypothesis on the geometry and only require that there is a
sufficiently rapid decay of local energy for the linear homogeneous
wave equation, which permits some domains that contain trapped rays.
The key step is to prove useful energy estimates involving the scaling
vector field for which the approach of the second author and Sogge
provides guidance.
\end{abstract}

\maketitle

\newsection{Introduction}
In this article, a lower bound of $c/\varepsilon^2$, where
$\varepsilon$ denotes the size of the Cauchy data, is established for
the lifespan of solutions to
quasilinear wave equations in exterior domains with Dirichlet boundary
conditions.  Here we examine nonlinearities which vanish to second
order and may depend on the solution $u$, rather than just its
derivatives, at the quadratic level.  The lifespan bound that is
established is an analog of that which was obtained in \cite{Lindblad}
in the absence of a boundary.  Similar lifespan bounds appeared in
\cite{DZ} (and could also be obtained via the methods of \cite{DMSZ})
exterior to star-shaped obstacles.  We permit much more general
geometries and only require that there is a sufficiently rapid decay
of local energy, with a possible loss of regularity.
Due, e.g., to the
seminal results \cite{Ikawa1, Ikawa2}, this includes some domains which contain
trapped rays.

Let us more specifically introduce the problem.  Let $\K\subset \R^3$
be a bounded domain with smooth boundary.  Note that we shall not
assume that $\K$ is connected.  We then examine the following
quasilinear wave equation exterior to $\K$
\begin{equation}
  \label{main}
  \begin{cases}
    \Box u(t,x)=Q(u,u',u''),\quad (t,x)\in \R_+\times\ext,\\
    u(t,\cd)|_{\bdy}=0,\\
    u(0,\cd)=f,\quad \partial_t u(0,\cd)=g.
  \end{cases}
\end{equation}
Here $\Box=\partial_t^2-\Delta$ is the d'Alembertian, and
$u'=\partial u = (\partial_t u, \nabla_x u)$ denotes the space-time gradient.  The
nonlinear term vanishes quadratically at the origin and is affine
linear in $u''$, thus yielding a quasilinear equation.  Without loss of generality, we may take
$0\in\K\subset\{|x|<1\}$, and we shall do this throughout.  While we shall
state the lifespan bound for the scalar equation \eqref{main}, the
methods shall fully permit, even multiple speed, systems.

The nonlinearity $Q$ can be expanded as
\[Q(u,u',u'')=A(u,u') +
B^{\alpha\beta}(u,u')\partial_\alpha\partial_\beta u\]
where $A(u,u')$ vanishes to second order at the origin and
$B^{\alpha\beta}$ are functions which are symmetric in $\alpha,\beta$
and vanish to the first
order at $(0,0)$.  Here we are using the summation convention where
repeated indices are implicitly summed from $0$ to $3$, where $x_0=t$
and $\partial_0=\partial_t$.

For simplicity of exposition, we shall truncate at the quadratic
level.  As we are examining a small amplitude problem, it is clear
that this shall not affect the long-time behavior.  Upon doing so,
we may write
\[Q(u,u',u'') = A(u,u') + b^{\alpha\beta}
u \partial_\alpha\partial_\beta u +
b^{\alpha\beta}_{\gamma} \partial_\gamma
u \partial_\alpha\partial_\beta u\]
where $b^{\alpha\beta}$ and $b^{\alpha\beta}_\gamma$ are constants
which are symmetric in $\alpha, \beta$ and $A(u,u')$ is a quadratic form.

In order to solve \eqref{main}, the Cauchy data are required to
satisfy certain compatibility conditions.  Letting
$J_ku=\{\partial_x^\alpha u\,:\,0\le |\alpha|\le k\}$, for a
formal $H^m$ solution $u$ of \eqref{main}, we can write $\partial_t^k
u(0,\cd)=\psi_k(J_kf, J_{k-1}g)$, $0\le k\le m$.  The functions
$\psi_k$, which depend on $Q$, are called compatibility functions.
The compatibility condition for $(f,g)\in H^m\times H^{m-1}$ simply
requires that $\psi_k$ vanishes on $\bdy$ for all $0\le k\le m-1$.
For smooth data, the compatibility conditions are said to hold to
infinite order if this condition holds for all $m$.  A more detailed
exposition on compatibility can be found in, e.g., \cite{KSS2}.

Our only assumption on the geometry of $\K$ is that the associated local energy
for solutions to the linear homogeneous wave equation
decays sufficiently rapidly.  For a clearer exposition, we shall
assume that there is an exponential decay of local energy with a loss
of a single degree of regularity, though it will be clear from the
proof that a sufficiently rapid polynomial decay with a fixed finite
loss of regularity will suffice.\footnote{To date, the authors are not aware
of examples of three dimensional domains for which there is
polynomial decay but not exponential decay, though the recent article
\cite{CW} provides the most compelling evidence to date that
such domains can be constructed.}  More specifically, we shall assume
that if $\Box u=0$ and if $\supp u(0,\cd), \partial_t u(0,\cd)\subset
\{|x|<10\}$, then
\begin{equation}\label{decay}
  \|u'(t,\cd)\|_{L^2(\ext\cap \{|x|<10\})}\lesssim e^{-a t}
  \sum_{|\alpha|\le 1} \|\partial^\alpha_x u'(0,\cd)\|_2
\end{equation}
for some $a>0$.  The notation $A\lesssim B$ indicates that there is a
positive unspecified constant $C$ (which may change from line to line)
so that $A\le CB$.  Moreover, this $C$ will implicitly be independent
of any parameters in our problem.

Local energy decay estimates such as \eqref{decay} have a long
history, which we shall only tersely describe.  For nontrapping
domains, one need not have the loss of regularity which appears in the
right.  See \cite{LMP} for star-shaped obstacles and \cite{MRS} for nontrapping
domains.  When there are trapped rays, it is known \cite{Ralston} that an
estimate such as \eqref{decay} cannot hold unless such a loss of
regularity is permitted.  Results in the positive direction in the
presence of trapped rays begin with \cite{Ikawa1, Ikawa2} where such estimates
are proved when $\K$ consists of multiple convex obstacles subject to
certain size/spacing conditions.  More recent results include
\cite{BZ}, \cite{C},
\cite{CdVP}, \cite{NZ},  \cite{WZ}.

We may now state our main theorem, which shows that for Cauchy data of
size $\varepsilon$ solutions to \eqref{main} must exist up to a
lifespan of $T_\varepsilon = c/\varepsilon^2$ for some small constant
$c$.
\begin{theorem}\label{thm1}
 Let $\K$ be a smooth, bounded set for which the exponential
 decay of local energy \eqref{decay} holds, and let $Q$ be as
 above.  Suppose that the Cauchy
 data $f,g\in C^\infty(\ext)$ are compactly supported and satisfy the compatibility conditions to
 infinite order.  Then there exist constants $N$ and $c$ so that if
 $\varepsilon$ is sufficiently small and
  \begin{equation}\label{data}
   \sum_{|\mu|\le N} \|\partial_x^\mu f\|_{2} + \sum_{|\mu|\le
     N-1}\|\partial_x^\mu g\|_2 \le \varepsilon,
  \end{equation}
then \eqref{main} has a unique solution $u\in
C^\infty([0,T_\varepsilon]\times \ext)$ where
\begin{equation}
  \label{lifespan}
T_\varepsilon = c/\varepsilon^2.
\end{equation}
\end{theorem}

For simplicity of exposition, we are assuming here that the Cauchy
data are compactly supported.  It is likely that it would suffice to
take the data to be small in certain weighted Sobolev norms.

On $\R_+\times\R^3$, this lifespan was first proved in
\cite{Lindblad}.  The dependence on the solution $u$ rather than just
its first and second derivatives inhibits the energy methods which are
typically employed to show such long time existence.  Indeed, when the
nonlinearity does not depend on $u$ at the quadratic level, solutions
are known to exist almost globally, i.e. with a lifespan of
$T_\varepsilon\approx \exp(c/\varepsilon)$.  See \cite{JK}.

The previous results \cite{KSS4, KSS, KSS2, KSS3}, \cite{MS3, MS,
  MS4}, \cite{MNS}, \cite{KK} have focused on proving
long-time existence for three dimensional wave equations in exterior domains when there
is no dependence on $u$, and \cite{MNS2} examines the case that there is
dependence at the cubic level and beyond.  Of particular note is the
paper \cite{MS3} where long-time existence results were first established
only assuming \eqref{decay}, and in particular, in domains which have
trapped rays.

The current direction of
research was initiated with the paper \cite{DZ} where the same lower
bound on the lifespan was shown exterior to star-shaped obstacles.  A
similar four dimensional problem was addressed in \cite{DMSZ},
and the techniques therein could also be applied to the exterior of
three dimensional star-shaped obstacles.  The current article extends
these results to much more general geometries.

The method of proof shall utilize Klainerman's method of invariant
vector fields \cite{Klainerman}.  This has been adapted to the exterior domain
setting with particularly notable contributions coming in \cite{KSS, KSS3}
 and \cite{MS3}.  To this end, we set 
\[Z=\{\partial_\alpha, \Omega_{ij}=x_i\partial_j-x_j\partial_i\,:\,
0\le \alpha\le 3,\, 1\le i<j\le 3\}\]
to be the generators of space-time translations and spatial
rotations.  We shall also denote $L=t\partial_t+r\partial_r$, which is
the scaling vector field.  A key fact is
that 
\[[\Box,Z]=0,\quad [\Box, L] = 2\Box.\]
These vector fields thus preserve $\Box u = 0$, in the sense that if
$u$ solves such an equation then $\Box Zu=0$ and $\Box Lu=0$.  The
Lorentz boosts $\Omega_{0k}=t\partial_k + x_k\partial_t$ have not yet
been mentioned despite playing a key role when the method is applied
on $\R_+\times \R^3$.  Though all of these vector fields have nice
commutation properties with $\Box$, only $\partial_t$ is guaranteed to
preserve the Dirichlet boundary conditions.  While the other members
of $Z$ do not preserve the boundary conditions, their coefficients are
bounded on the compact obstacle, and thus, approximately do.  Indeed,
these can be handled using elliptic regularity arguments and localized
energy estimates as was initiated in \cite{KSS}.  On the other hand,
the Lorentz boosts have unbounded normal component on $\bdy$ and seem
inadmissible for such boundary value problems.  It is also worth
noting that the Lorentz boosts have an associated speed, and they only
commute with the d'Alembertian of the same speed, which renders them
also difficult to use for multiple speed systems.  While the scaling
vector field has a bounded normal component on $\bdy$, for long-time
problems, its coefficients can be large in any neighborhood of the
boundary.  For this reason, we shall be required to use relatively few
$L$ compared to the vector fields $Z$.  This is an idea which
originated in \cite{KSS3} and is further displayed in \cite{MS3},
\cite{MNS, MNS2}. 

The star-shaped assumption on the geometry of $\K$ arises in \cite{DZ}
in order to prove energy estimates involving the scaling vector field
$L$.  Indeed, using ideas akin to those from \cite{KSS3} which is in
turn reminiscent of \cite{Morawetz}, it is shown that the worst
boundary term (in terms of $t$ dependence) 
which arises when proving an energy inequality for $Lu$
has a beneficial sign.  Developing an alternative method for handling
these boundary terms was one of the major innovations of \cite{MS3},
and this article largely represents a combination of ideas from
\cite{DZ} and \cite{MS3}.

The star-shaped assumption arises in \cite{DMSZ} in a related, though
different, way.  Long-time existence is shown there using only the
vector fields $Z$.  This is accomplished by employing a class of localized
energy estimates which are known to hold for small perturbations of
the d'Alembertian exterior to star-shaped obstacles \cite{MS}, \cite{MT2} and
iterating in a fashion which is akin to \cite{KSS} and \cite{MS}.  See
also \cite{MS4} for a further example of how a star-shaped assumption and
the broader class of available localized energy estimates can simplify
arguments.

The remainder of the article is organized as follows.  In Section 2,
we gather our main energy and localized energy estimates.  These
largely represent a combination of the main estimates used in
\cite{DZ} as well as the techniques developed in \cite{MS3} to permit
the use of the scaling vector field when the obstacle is not
star-shaped.  In Section 3, we establish the main decay estimates
which we shall utilize.  The principal piece here is a $L^1-L^\infty$
estimate which is akin to those of H\"ormander
\cite{HormanderL1Linfty} as was adapted to remove the dependence on
the Lorentz boosts by \cite{KSS3}.  In Section 4, we prove the long
time existence given by Theorem \ref{thm1}.

\newsection{Energy estimates}
In this section, we shall gather the main $L^2$ estimates which we
shall require.  These will primarily be energy estimates as well as
weighted $L^2_tL^2_x$ localized energy estimates for the solution and
vector fields applied to the solution.  That is, these will
be variants of the energy estimate
and the localized energy estimate
\begin{multline}\label{le}\sup_{t\in [0,T]} \|u'(t,\cd)\|_2 +
  \sup_{R\ge 1} R^{-1/2}
  \|u'\|_{L^2_{t,x}([0,T]\times \{|x|<R\})}\lesssim \|u'(0,\cd)\|_2 \\+
  \inf_{\Box u = f+g}\Bigl(
\int_0^T \|f(s,\cd)\|_2\,ds + \sum_{j\ge 0} \|\la x\ra^{1/2}
g\|_{L^2_{t,x}([0,T]\times \{\la x\ra \approx 2^j\})}\Bigr)\end{multline}
which are known to hold on $\R_+\times \R^3$.  Estimates of this
latter form originated in \cite{Morawetz2} and have subsequently
appeared in, e.g., \cite{HY}, \cite{KSS, KSS3}, \cite{KPV}, \cite{SS}, \cite{Sterb},
\cite{Strauss}.  Their particular utility for proving long-time existence in exterior
domains was first recognized in \cite{KSS}, and they have played a primary
role in nearly every such proof since.  Also, see, e.g., \cite{MS,
  MS4} and \cite{MT1, MT2}.

\subsubsection{Estimates for $\|u\|_{L^2_x}$ on $\R_+\times\R^3$}\label{bdylessSection}
Here we shall gather the boundaryless $L^2$ estimates on $u$, rather
than on $u'$, which we shall utilize in the sequel.  We shall only
require these in the boundaryless case as the Dirichlet boundary
conditions permit the control 
\begin{equation}\label{locControl}\|u\|_{L^2_x(\{x\in
  \ext\,:\,|x|<2\})} \lesssim \|u'\|_{L^2_x(\{x\in
  \ext\,:\,|x|<2\})}\end{equation}
and when a cutoff which vanishes on $\{|x|<1\}$
is applied to $u$, then it suffices to examine a boundaryless
equation.  The majority of these estimates were also utilized in
\cite{DZ}.  In the sequel, we shall abbreviate $L^2_x(\{x\in
\ext\,:\,|x|<2\})$ as $L^2_x(|x|<2)$.

We first state the variant of the localized energy estimates which we
shall employ.

\begin{lemma} 
  Let $u$ be a smooth function which vanishes for large $|x|$ at each
  time $t$.  Then for $T\ge 1$
  \begin{multline}\label{DZle}
\|\la x\ra^{-3/4} u'\|_{L^2_{t,x}([0,T]\times\R^3)}+ T^{-1/4} \|\la x\ra ^{-1/4} u'\|_{L^2_{t,x}([0,T]\times
      \R^3)}
\\\lesssim \|u'(0,\cd)\|_2 
+
  \inf_{\Box u = f+g}\Bigl(
\int_0^T \|f(s,\cd)\|_2\,ds + \sum_{j\ge 0} \|\la x\ra^{1/2}
g\|_{L^2_{t,x}([0,T]\times \{\la x\ra \approx 2^j\})}\Bigr).
%
  \end{multline}
\end{lemma}

To obtain \eqref{DZle}, we first note that over $|x|>T$ it follows trivially from
the energy inequality.
To finish the proof, one need only dyadically decompose $|x|<T$ and
apply \eqref{le}.  The interested reader can find an alternate proof
in \cite{DZ}.

We shall then use the following weighted Sobolev-type estimate from
\cite{DZ}, \cite{DMSZ}.
\begin{lemma}
For $n\ge 3$ and $h\in C_0^\infty(\R^n)$, 
\[\|h\|_{\dot{H}^{-1}(\R^n)}\lesssim \|h\|_{L^{\frac{2n}{n+2}}(|x|<1)}
+ \||x|^{-\frac{n-2}{2}} h\|_{L^1_rL^2_\omega(|x|>1)}.\]
\end{lemma}

Here and throughout, the mixed norm represents
\[\|f\|_{L^p_r L^q_\omega(\R^n)} = \Bigl(\int_0^\infty
\Bigr[\int_{\S^{n-1}} |f(r\omega)|^q\,d\omega\Bigr]^{p/q}\, r^{n-1}\,dr\Bigr)^{1/p}.\]

By applying the energy inequality and \eqref{DZle} to the Riesz
transforms of the solution and subsequently applying the preceding
lemma, we obtain:

\begin{proposition}
  Let $u$ be a smooth function which vanishes for large $|x|$ at each
  time $t$.  Then for $T\ge 1$
  \begin{multline}\label{noD}
    \|u\|_{L^\infty_tL^2_x([0,T]\times \R^3)} + T^{-1/4} \|\la x\ra^{-1/4}
    u\|_{L^2_{t,x}([0,T]\times \R^3)} \lesssim \|u(0,\cd)\|_2 +
    \|\partial_tu(0,\cd)\|_{\dot{H}^{-1}}
\\ + \int_0^T \||x|^{-1/2} \Box u(s,\cd)\|_{L^1_rL^2_\omega(|x|>1)}\,ds +
\int_0^T \|\Box u(s,\cd)\|_{L^{6/5}(|x|<1)}\,ds.
  \end{multline}
\end{proposition}

Further details of the proof can again be found in \cite{DZ} and
\cite{DMSZ}.

The following proposition will be used to control commutator terms
when the solution in the exterior domain is cutoff away from the
obstacle.  It appears implicitly in \cite[Section 4]{DZ} and utilizes
arguments akin to those which appeared in \cite{SS}, \cite{KSS} which
rely on Huygens' principle.  In higher dimensions, an alternate proof
which does not rely on Huygens' principle appeared in \cite{DMSZ}.

\begin{proposition}
  Let $u$ be a smooth solution to $\Box u = G$, $u=0$ for $t\le 0$.
  Suppose further that $G(s,x)=0$ unless $|x|<3$.  Then,
  \begin{equation}
    \label{noDcomm}
    \|u\|_{L^\infty_tL^2_x([0,T]\times \R^3)} + T^{-1/4} \|\la x\ra^{-1/4}
    u\|_{L^2_{t,x}([0,T]\times \R^3)} \lesssim
    \|G\|_{L^2_{t,x}([0,T]\times \R^3)}.
  \end{equation}
\end{proposition}

In one case, we shall need an improvement on the estimate \eqref{noD} when
the forcing term is in divergence form.  This follows easily from
ideas of \cite{Hormander}, \cite{Lindblad}.  See also \cite{MS5} for
an application in a context similar to the current study.

\begin{proposition}\label{propDivForm}
Let $v$ be a smooth solution to
\[\begin{cases}
  \Box v = \sum_0^3 a_j\partial_j G,\quad (t,x)\in \R_+\times \R^3,\\
v(0,\cd)=\partial_t v(0,\cd)=0.
\end{cases}\]
Then,
\begin{equation}
  \label{divForm}
  \sup_{t\in [0,T]}\|v(t,\cd)\|_2 + \la T\ra^{-1/4} \|\la x\ra^{-1/4}
  v\|_{L^2_{t,x}([0,T]\times \R^3)} \lesssim
  \|G(0,\cd)\|_{\dot{H}^{-1}} + \int_0^T \|G(s,\cd)\|_2\,ds.
\end{equation}
\end{proposition}

\subsubsection{Energy estimates on $\R_+\times\ext$}
In this section, we examine the fixed time, $L^2$ energy estimates
which will be used in the sequel.  As we are proving long-time
existence for quasilinear equations, we shall require such estimates
for small perturbations of the d'Alembertian.  Such an estimate is
well-known for solutions satisfying Dirichlet boundary conditions.
When vector fields are, however, applied to the solution and these
boundary terms no longer vanish, some extra care is required,
particularly when the time dependent vector field $L$ occurs.  In the
remainder of this section, we are merely gathering results from
\cite[Section 2]{MS3}, and the interested reader is referred therein
for detailed proofs.

In particular, we shall be studying smooth solutions $u$ to
\begin{equation}\label{pert}
  \begin{cases}
    \Box_\gamma u = F,\quad (t,x)\in \R_+\times\ext,\\
    u|_{\bdy}(t,\cd)=0,\\
    u(0,\cd)=f,\quad \partial_t u(0,\cd)=g
  \end{cases}
\end{equation}
where
\[\Box_\gamma u = (\partial_t^2 -\Delta)u +
\gamma^{\alpha\beta}(t,x)\partial_\alpha\partial_\beta u.\]
The perturbation is taken to satisfy
$\gamma^{\alpha\beta}=\gamma^{\beta\alpha}$ as well as 
\begin{equation}\label{gammaDecay1}
  \|\gamma^{\alpha\beta}(t,\cd)\|_\infty \le \frac{\delta}{1+t},\quad
  0<\delta\ll 1.
\end{equation}
We set $e_0(u)$ to be the energy form
\[
  e_0(u)=|u'|^2 + 2\gamma^{0\alpha} \partial_t u \partial_\alpha u -
  \gamma^{\alpha\beta}\partial_\alpha u\partial_\beta u.
\]

Our first estimate concerns
\[E_M(t)=E_M(u)(t)=\int_{\ext} \sum_{j=0}^M e_0(\partial_t^j
u)(t,x)\,dx.\]
As $\partial_t^j$ preserves the Dirichlet boundary conditions,
standard energy methods yield
\begin{lemma}\label{lemma2.1}
  Fix $M=0,1,2,\dots$ and assume that the $\gamma^{\alpha\beta}$ are
  as above.  Suppose that $u\in C^\infty$ solves \eqref{pert} and
  vanishes for large $|x|$ for every $t$.  Then 
  \begin{equation}
    \label{dtEnergy}
  \partial_t E_M^{1/2}(t)\lesssim \sum_{j=0}^M \|\Box_\gamma
  \partial_t^j u(t,\cd)\|_2 + \|\gamma'(t,\cd)\|_\infty E_M^{1/2}(t).
  \end{equation}
\end{lemma}
Here, we have used the notation
\[ \|\gamma'(t,\cd)\|_\infty = \sum_{\alpha,\beta,\mu=0}^3
\|\partial_\mu \gamma^{\alpha\beta}(t,\cd)\|_\infty.\]
In the sequel, we shall frequently use the fact that
\[\int_{\ext} e_0(v)(t,x)\,dx \approx \|v'(t,\cd)\|^2_2\]
if \eqref{gammaDecay1} holds with $\delta$ sufficiently small.

From these $L^2$ estimates for $\partial_t$ applied to the solution
$u$, estimates for $\partial_x^\mu u$ shall be obtained via elliptic
regularity.  The key lemma is
\begin{lemma}\label{lemma2.3}
  For $j, N=0,1,2,\dots$ fixed and for $u\in C^\infty(\R_+\times\ext)$
  solving \eqref{pert} and vanishing for large $|x|$ for each $t$, we
  have
  \begin{equation}\label{ellReg}
    \sum_{|\mu|\le N} \|L^j \partial^\mu u'(t,\cd)\|_2 \lesssim
    \sum_{\substack{k+l\le j+N\\ l\le j}} \|L^l\partial_t^k
    u'(t,\cd)\|_2 + \sum_{\substack{|\mu|+l\le N+j-1\\l\le j}}
    \|L^l\partial^\mu \Box u(t,\cd)\|_2.
  \end{equation}
\end{lemma}

In order to prove estimates involving $L$, we set
$\tilde{L}=\eta(x)r\partial_r + t\partial_t$ where $\eta\in
C^\infty(\R^3)$ with $\eta(x)=0$ for $x\in \K$ and $\eta(x)=1$ for
$|x|>1$.  We note that $\tilde{L}$ now preserves the Dirichlet
boundary conditions.  It, however, no longer commutes with $\Box$.  The
commutators will be controlled using a combination of \eqref{decay}
and Huygens' principle for the associated boundaryless d'Alembertian,
which will be stated later.  We set
\[X_{k,j} = \int e_0(\tilde{L}^k \partial_t^j u)(t,x)\,dx.\]
Associated to this energy, we have the following estimate.

\begin{proposition}
  Let $u\in C^\infty$ solve \eqref{pert} with $\gamma^{\alpha\beta}$
  as above and vanish for large $|x|$
  for each $t$.  Then,
  \begin{multline}
    \label{X}
    \partial_t X^{1/2}_{k,j} \lesssim \|\gamma'(t,\cd)\|_\infty X^{1/2}_{k,j}
    +\|\tilde{L}^k \partial_t^j
    \Box_\gamma u(t,\cd)\|_2 + \|[\tilde{L}^k\partial_t^j,
    \gamma^{\alpha\beta}\partial_\alpha\partial_\beta]u(t,\cd)\|_2 \\+
    \sum_{l\le k-1} \|L^l\partial_t^j \Box u(t,\cd)\|_2 +
    \sum_{\substack{l+|\mu|\le j+k\\l\le k-1}} \|L^l \partial^\mu u'(t,\cd)\|_{L^2(|x|<1)}.
  \end{multline}
\end{proposition}

In the sequel, we shall be choosing $\gamma$ so that
\begin{equation}\label{gammaDecay2}\|\gamma'(t,\cd)\|_\infty \le
  \frac{\delta}{1+t}.\end{equation}
By doing so, it will be easy to bootstrap the term involving
$\|\gamma'\|_\infty$ upon integration over $[0,T]$ when $T$ is
appropriately bounded in terms of $\delta$.


We finally state an energy estimate involving the full set of vector
fields.  Here, the associated boundary terms involve a loss of
regularity, but they no longer involve the rotations.  These boundary
terms will be controlled using localized energy estimates which follow.

\begin{proposition}
 For fixed $N_0$ and $m_0$, set
\[Y_{N_0,m_0}(t) = \sum_{\substack{|\mu|+k\le N_0+m_0\\k\le m_0\\|\nu|=1}} \int
e_0(L^k Z^\mu \partial^\nu u)(t,x)\,dx.\]
Suppose that \eqref{gammaDecay1} and \eqref{gammaDecay2} hold for
$\delta$ sufficiently small.  Then
\begin{multline}
  \label{LZenergy}
\partial_t Y_{N_0,m_0}\lesssim
Y_{N_0,m_0}^{1/2}\sum_{\substack{|\mu|+k\le N_0+m_0\\k\le
    m_0\\|\nu|=1}}\|\Box_\gamma L^k Z^\mu \partial^\nu u(t,\cd)\|_2 \\+
\|\gamma'(t,\cd)\|_\infty Y_{N_0,m_0} + \sum_{\substack{|\mu|+k\le
    N_0+m_0+2\\k\le m_0}} \|L^k \partial^\mu u'(t,\cd)\|^2_{L^2(|x|<1)}.
\end{multline}
\end{proposition}

The above proposition contains a slight modification of what appeared
previously in \cite{MS3}.  There one did not need to distinguish
between the vector fields $Z$ and the one derivative
$\partial$ in the definition of $Y_{N_0,m_0}$.  Here we need
this slight bit of additional precision.  The proof follows the same
argument.  One needs to only apply standard energy methods to the
principle terms and a trace theorem to the boundary terms which result
upon doing such integrations by parts.

\subsubsection{Localized energy estimates and boundary term estimates on $\R_+\times\ext$}
In this section, we collect two additional results from \cite{MS3}.
The reader is referred there for proofs.
Both estimates will concern solutions to the Dirichlet-wave equation
  \begin{equation}\label{inhomwave}
    \begin{cases}
      \Box w = G(t,x),\quad (t,x)\in \R_+\times\ext,\\
      w|_{\bdy}(t,\cd)=0,\\
      w(t,x)=0,\quad t\le 0.
    \end{cases}
  \end{equation}

When the estimates of Section \ref{bdylessSection} are applied away
from the obstacle, the following, which strongly depends on
\eqref{decay}, is used to handle the behavior near the boundary.  This
is also used to control the boundary term that appears in
\eqref{LZenergy}.  For
notational convenience, we set $S_T=[0,T]\times \ext$.

\begin{proposition}  Suppose that $\K\subset \{|x|<1\}$ satisfies
  \eqref{decay}, and suppose that
$w\in C^\infty$ solves \eqref{inhomwave}.  Then, for fixed $N_0$ and
$m_0$, if $w$ vanishes for large $|x|$ for every fixed $t$, 
\begin{multline}
  \label{locEnergyExt}
\sum_{\substack{|\mu|+j\le N_0+m_0\\j\le m_0}} \|L^j \partial^\mu
w'\|_{L^2_{t,x}(S_T\cap \{|x|<10\})} \lesssim \int_0^T
\sum_{\substack{|\mu|+j\le N_0+m_0+1\\j\le m_0}} \|\Box
L^j \partial^\mu w(s,\cd)\|_2\,ds
\\+ \sum_{\substack{|\mu|+j\le N_0+m_0-1\\j\le m_0}} \|\Box
L^j \partial^\mu w\|_{L^2_{t,x}(S_T)}.
\end{multline}
\end{proposition}

The second estimate shall be used to control the boundary term which
arises in \eqref{X}.  The contributions from behavior near the
obstacle are controlled using \eqref{decay} and that away from the
boundary follows from sharp Huygens' principle after passing to a
properly related boundaryless equation.  See \cite{MS3}.
\begin{proposition}
  Let $\K\subset \{|x|<1\}$ satisfy \eqref{decay}, and suppose that
  $w\in C^\infty$ solves \eqref{inhomwave}.  Then for fixed $N_0$ and
  $m_0$ and $t>2$,
  \begin{multline}
    \label{ms_bdy}
    \sum_{\substack{|\mu|+j\le N_0+m_0\\ j\le m_0}} \int_0^t
    \|L^j \partial^\mu w'(s,\cd)\|_{L^2(|x|<2)}\,ds \\\lesssim
    \sum_{\substack{|\mu|+j\le N_0+m_0+1\\j\le m_0}} \int_0^t
    \Bigl(\int_0^s \|L^j \partial^\mu \Box w(\tau,\cd)\|_{L^2(||x|-(s-\tau)|<10)}\,d\tau\Bigr)\,ds.
  \end{multline}
\end{proposition}


\newsection{Pointwise decay estimates}
In this section, we gather the main decay estimates which will permit
the necessary integrability to gain long-time existence.  The first is
a variant on standard weighted $L^1$-$L^\infty$ estimates (see,
e.g., \cite{Htext}, \cite{So2}).

\begin{proposition}\label{l1linfProp}
  Suppose that $w$ is a solution to the scalar inhomogeneous wave
  equation \eqref{inhomwave},
and suppose that $\K\subset \{|x|<1\}$ is so that the decay of local energy
\eqref{decay} holds.  Then
\begin{multline}\label{l1linf}
  (1+t+|x|) |Z^\mu w(t,x)|\lesssim \int_0^t\int_\ext
  \sum_{\substack{|\nu|+k\le |\mu|+7\\k\le 1}} |L^k
  Z^\nu G(s,y)|\,\frac{dy\,ds}{|y|}
\\+\int_0^t \sum_{\substack{|\nu|+k\le |\mu|+4\\k\le 1}}
\|L^k \partial^\nu G(s,\cd)\|_{L^2(|x|<2)}\,ds.
\end{multline}
\end{proposition}

This is essentially \cite[Theorem 4.1]{KSS3}, though there solutions
are studied exterior to star-shaped obstacles.  For such domains, the
associated version of \eqref{decay} does not require a loss of
regularity.  In the current setting, as we allow for the possibility
of trapped rays in our exterior domain and as such \eqref{decay} has
an associated loss of regularity, the right hand side of \eqref{l1linf} reflects
an additional vector field.  See, also, \cite[Theorem 3.1]{MS3}.

The second decay estimate is a weighted Sobolev lemma.  See
\cite{Klainerman2}.

\begin{lemma}\label{wtdSoblemma}
  Suppose that $h\in C^\infty(\R^3)$.  Then for $R\ge 1$, 
  \begin{equation}
    \label{wtdSob}
   \|h\|_{L^\infty(R/2<|x|<R)}\lesssim R^{-1}\sum_{|\alpha|\le 2}
   \|Z^\alpha h\|_{L^2(R/4<|x|<2R)},
  \end{equation}
and
 \begin{equation}
    \label{wtdSob2}
   \|h\|_{L^\infty(R<|x|<R+1)}\lesssim R^{-1}\sum_{|\alpha|\le 2}
   \|Z^\alpha h\|_{L^2(R-1<|x|<R+2)}.
  \end{equation}
\end{lemma}

After localizing to the annulus, these estimates follow simply by
applying Sobolev embedding on $\R\times \S^2$ and then adjusting the
volume element to match that of $\R^3$.

\newsection{Proof of Theorem \ref{thm1}}
Here we prove Theorem~\ref{thm1} via iteration.  We shall first use
local existence theory to reduce to the case of vanishing initial data.

Indeed by invoking, e.g., the local existence theory established in
\cite{KSS2}, if $\varepsilon$ in \eqref{data} is sufficiently small
and $N$ is sufficiently large,
then the existence of a smooth solution
for $t\in [0,2]$ satisfying
\begin{equation}
  \label{local}
  \sup_{0\le t\le 2} \sum_{|\mu|\le 102} \|\partial^\mu
  u(t,\cd)\|_{2}\le C\varepsilon
\end{equation}
is guaranteed.

We shall use this local solution to reduce to the case of vanishing
initial data.  To this end, let $\eta\in C^\infty(\R)$ with
$\eta(t)\equiv 1$ for $t\le 1/2$ and $\eta(t)\equiv 0$ for $t>1$.
Then $u_0=\eta u$ solves
\[\Box u_0 = \eta Q(u,u',u'') + [\Box,\eta]u.\]
And solving \eqref{main} is then equivalent to showing that $w=u-u_0$
solves
\begin{equation}\label{reduced}
  \begin{cases}
    \Box w = (1-\eta)Q(u_0+w,(u_0+w)',(u_0+w)'') -
         [\Box,\eta](u_0+w),\\
    w|_\bdy = 0,\\
    w(0,x)=\partial_tw(0,x)=0.
  \end{cases}
\end{equation}

We solve \eqref{reduced} via iteration.  In particular, we let
$w_0\equiv 0$, and recursively define $w_m$ to solve
\[
\begin{cases}
  \Box w_m = (1-\eta) Q(u_0+w_{m-1},(u_0+w_{m-1})', (u_0+w_m)'') - [\Box,\eta]u,\\
w_m|_{\bdy}=0,\\
w_m(0,x)= \partial_t w_m(0,x)=0.
\end{cases}
\]

Our first goal is to show a form of boundedness.
We set
\begin{multline}\label{M}
M_m(T)=\sup_{t\in [0,T]}\sum_{|\mu|\le 100}\|\partial^\mu
w_m'(t,\cd)\|_2 
+ \sum_{|\mu|\le 95} \|\la x
\ra^{-3/4}\partial^\mu w_m'\|_{L^2_{t,x}(S_T)}\\
+\sup_{t\in [0,T]} \sum_{\substack{|\mu|\le 90\\|\nu|\le 2}}
\|Z^\mu\partial^\nu w_m(t,\cd)\|_2 
+ \sum_{\substack{|\mu|\le
    90\\|\nu|\le 1}} \la T\ra^{-1/4} \|\la x\ra^{-1/4}
Z^\mu\partial^\nu w_m\|_{L^2_{t,x}(S_T)}\\
+\sup_{0\le t\le T} \sum_{|\mu|\le 80} \|L\partial^\mu w'_m(t,\cd)\|_2
+ \sum_{|\mu|\le 75} \|\la x\ra^{-3/4} L\partial^\mu
w_m'\|_{L^2_{t,x}(S_T)}
\\+ \sup_{0\le t\le T} \sum_{\substack{|\mu|\le 70\\|\nu|\le 2}}
  \|LZ^\mu \partial^\nu w_m(t,\cd)\|_2
+ \la T\ra^{-1/4} \sum_{\substack{|\mu|\le 70\\|\nu|\le 1}} \|\la
  x\ra^{-1/4} LZ^\mu \partial^\nu w_m\|_{L^2_{t,x}(S_T)}
\\+ \sup_{0\le t\le T} (1+t) \sum_{|\mu|\le 60} \|Z^\mu w_m(t,\cd)\|_\infty.
\end{multline}
If $M_0(T)$ denotes the above quantity with $w_m$ replaced by
$u_0$, then \eqref{local} and finite propagation speed guarantees the existence of a constant $C_0$
so that
\begin{equation}\label{basecase}M_0(T)\le
  C_0\varepsilon.\end{equation}
We wish to inductively show that there is a uniform constant $C_1$ so that
\begin{equation}\label{bddness}M_m(T)\le C_1\varepsilon.\end{equation}


We label the terms of $M_m(T)$ by $I, II, \dots, IX$.

{\bf \em Bound for $I$:}
Here we shall apply \eqref{dtEnergy} and \eqref{ellReg} ($j=0$) with
\begin{equation}\label{gamma}\gamma^{\alpha\beta} = - (1-\eta)\Bigl[b^{\alpha\beta} (u_0+w_{m-1}) +
b^{\alpha\beta}_\sigma \partial_\sigma(u_0+w_{m-1})\Bigr].\end{equation}
By the inductive hypothesis as well as \eqref{local}, we
have \eqref{gammaDecay1} and \eqref{gammaDecay2} with $\delta =
C_1\varepsilon$.  Upon integrating \eqref{dtEnergy}, applying
\eqref{lifespan}, bootstrapping, and utilizing \eqref{ellReg}, it
suffices to bound
\[\sum_{j\le 100} \int_0^T \|\Box_\gamma \partial_t^j w_m(t,\cd)\|_2\,dt
+ \sup_{t\in [0,T]}\sum_{|\mu|\le 99} \|\partial^\mu \Box w_m(t,\cd)\|_2.\]
Here we note that
\begin{multline*}
  \sum_{j\le 100} |\Box_\gamma \partial^j_t w_m|
\lesssim \sum_{|\nu|\le
    50} |\partial^\nu (u_0+w_{m-1})| \Bigl[
  \sum_{|\mu|\le 100} |\partial^\mu (u_0+w_m)'|
+ \sum_{|\mu|\le 102}
|\partial^\mu u_0|\Bigr]
\\ 
+ \sum_{|\nu|\le
    51} |\partial^\nu (u_0+w_{m})'|
  \sum_{|\mu|\le 100} |\partial^\mu (u_0+w_{m-1})'|
+ \sum_{|\nu|\le
    51} |\partial^\nu (u_0+w_{m-1})|
  \sum_{|\mu|\le 100} |\partial^\mu (u_0+w_{m-1})'|
\\+ |u_0+w_{m-1}|^2
+ \sum_{|\mu|\le 100} |\partial^\mu [\Box,\eta]u|.
\end{multline*}
By using term $I$, $III$, and $IX$ of \eqref{M} as well as \eqref{local}, we have
\begin{multline*}
  \sum_{j\le 100} \int_0^T \|\Box_\gamma \partial^j_t w_m(t,\cd)\|_2\,
  dt\lesssim
 (M_0(T)+M_{m-1}(T))(M_0(T)+M_m(T))  \int_0^T (1+t)^{-1}\,dt \\
 + (M_0(T)+M_{m-1}(T))^2 \int_0^T (1+t)^{-1}\,dt 
+\varepsilon.
\end{multline*}
We can argue similarly to control 
\[\sup_{t\in [0,T]}\sum_{|\mu|\le 99} \|\partial^\mu \Box
w_m(t,\cd)\|_2.\]
By doing so, it follows that
\[I\le
C(M_0(T)+M_{m-1}(T))(M_0(T)+M_{m-1}(T)+M_m(T))\log(2+T) +C_2\varepsilon\]
provided that $C_2$ is a constant which is chosen sufficiently large relative to the
constant in \eqref{local}.  At this point, we shall permit $C_2$ to
change from line to line but note that $C_2$ is completely independent
of important parameters such as $m$, $C_1$, $\varepsilon$, and $T$.


{\bf\em Bound for $II$:}
Here we fix a smooth cutoff $\beta$ which is identically $1$ on
$|x|<2$ and vanishes for $|x|>3$.  For the multi-index $\mu$ fixed, we
first examine $(1-\beta)\partial^\mu w_m$,
which solves the boundaryless wave equation
\[\Box (1-\beta)\partial^\mu w_m = (1-\beta)\partial^\mu \Box w_m -
[\Box,\beta] \partial^\mu w_m\]
with vanishing initial data.  To this, we apply \eqref{DZle}.  Thus,
in order to control $II$, 
we see that it suffices to bound
\[\sum_{|\mu|\le 95} \int_0^T \|\partial^\mu \Box w_m(s,\cd)\|_2\,ds +
\sum_{|\mu|\le 95} \|\partial^\mu w_m'\|_{L^2_{t,x}(S_T\cap
  \{|x|<3\})}.\]
Here we have applied \eqref{locControl} to the lower order piece of
the commutator.
To control the latter term, we utilize \eqref{locEnergyExt}, which
reduces the bound for $II$ to controlling
\begin{equation}\label{II_RHS}\sum_{|\mu|\le 96} \int_0^T \|\partial^\mu \Box w_m(s,\cd)\|_2\,ds +
\sum_{|\mu|\le 94} \|\partial^\mu \Box
w_m\|_{L^2_{t,x}(S_T)}.\end{equation}
As
\begin{multline}
  \label{productrule}
\sum_{|\mu|\le 96} |\partial^\mu \Box w_m| \lesssim
\sum_{|\nu|\le
    49} |\partial^\nu (u_0+w_{m-1})|
  \sum_{|\mu|\le 97} |\partial^\mu (u_0+w_m)'|\\ 
+ \sum_{|\nu|\le
    49} |\partial^\nu (u_0+w_{m})'|
  \sum_{|\mu|\le 96} |\partial^\mu (u_0+w_{m-1})'|
+ \sum_{|\nu|\le
    49} |\partial^\nu (u_0+w_{m-1})|
  \sum_{|\mu|\le 96} |\partial^\mu (u_0+w_{m-1})'|
\\+ |u_0+w_{m-1}|^2
+ \sum_{|\mu|\le 96} |\partial^\mu [\Box,\eta]u|,
\end{multline}
it follows from \eqref{M} ($I$, $III$, and $IX$) and \eqref{local} that
\begin{multline*}
\sum_{|\mu|\le 97} \int_0^T \|\partial^\mu \Box w_m(s,\cd)\|_2\,ds
\lesssim (M_0(T)+M_{m-1}(T))(M_0(T)+M_m(T)) \int_0^T s^{-1}\,ds 
\\+ (M_0(T)+M_{m-1}(T))^2\int_0^T s^{-1}\,ds + \varepsilon.
\end{multline*}
The second term in \eqref{II_RHS} is simpler and can be controlled
similarly.  It follows that
\[II\le C(M_0(T)+M_{m-1}(T))(M_0(T)+M_{m-1}(T)+M_m(T))\log(2+T) +C_2\varepsilon.\]

{\bf \em Bounds for $III$ and $IV$, $|\nu|=0$:}
The primary estimate which shall be utilized is \eqref{noD}.  This
meshes well with every nonlinear term with the exception of those
involving second derivatives, which is more difficult as we cannot
properly control the second derivatives in the weighted $L^2_{t,x}$
spaces.  To get around this, we write the worst terms in divergence
form and utilize \eqref{divForm}.

To do
so, we fix $\beta$ as above.  Over $|x|<3$, due to the Dirichlet
boundary conditions, we have that $\|w_m(s,\cd)\|_{L^2(|x|<3)} \lesssim
\|w_m'(s,\cd)\|_{L^2(|x|<3)}$.  Moreover, on such a compact set, the
coefficients of $Z$ are bounded, and
$\|Zw_m(s,\cd)\|_{L^2(|x|<3)}\lesssim \|w_m'(s,\cd)\|_{L^2(|x|<3)}$.   
Such terms are subject to the bound established for $II$ above.
Thus, it will suffice to control $(1-\beta)Z^\mu
w_m$ in the appropriate norms.

We note that for $\mu$ fixed, $(1-\beta)Z^\mu w_m$ solves the
boundaryless equation
\[\Box (1-\beta)Z^\mu w_m = (1-\beta)Z^\mu \Box w_m -[\Box,\beta]Z^\mu
w_m\]
and that the latter term is supported in $\{|x|<3\}$.  We must further
decompose the first term in the right.  We write
\begin{multline*}
(1-\beta)Z^\mu \Box w_m = \partial_\alpha\Bigl[(1-\eta)(1-\beta) \Bigl(b^{\alpha\beta}
(u_0+w_{m-1}) \partial_\beta Z^\mu (u_0+w_m) \\+
b^{\alpha\beta}_\gamma \partial_\gamma (u_0+w_{m-1}) \partial_\beta Z^\mu (u_0+w_m)\Bigr)\Bigr] 
+G_m(t,x).
\end{multline*}
The key here is that the former term falls into the context of
Proposition \ref{propDivForm} while the latter term does not contain
the case where the full number of vector fields lands on the term
containing second derivatives.

By applying \eqref{divForm}, \eqref{noD} and \eqref{noDcomm}, as well as the comments in
the preceding paragraph, we obtain
\begin{multline}\label{zeroCase}
  \sup_{t\in [0,T]} \sum_{|\mu|\le 90} \|(1-\beta) Z^\mu
  w_m(t,\cd)\|_2 + \sum_{|\mu|\le 90} \la T\ra^{-1/4} \|\la
  x\ra^{-1/4} (1-\beta)Z^\mu w_m\|_{L^2_{t,x}(S_T)}
\\\lesssim \int_0^T
\|\la x\ra^{-1/2}G_m(s,\cd)\|_{L^1_rL^2_\omega}\,ds 
+\sum_{|\mu|\le 90}\sum_{|\theta|\le 1}\int_0^T \|
|\partial^\theta (u_0+w_{m-1})| |Z^\mu (u_0+w_m)'|\|_2\,ds
\\+ \sum_{|\mu|\le 90} \|\partial^\mu w_m'\|_{L^2_{t,x}([0,T]\times\{|x|<3\})}.
\end{multline}
The last term is controlled by $II$, which was appropriately bounded
in the previous section.  

We first examine the first term in the right.  Analogous to
\eqref{productrule}, we have
\begin{multline}
  \label{productrule2}
|G_m| \lesssim
\sum_{|\nu|\le
    46} |Z^\nu (u_0+w_{m-1})|
  \sum_{|\mu|\le 90} |Z^\mu (u_0+w_m)'|\\ 
+ \sum_{|\nu|\le
    46} |Z^\nu (u_0+w_{m})'|
  \sum_{\substack{|\mu|\le 90\\|\theta|\le 1}} |Z^\mu \partial^\theta (u_0+w_{m-1})|
\\+ \sum_{|\nu|\le
    46} |Z^\nu (u_0+w_{m-1})|
  \sum_{\substack{|\mu|\le 90\\|\theta|\le 1}} |Z^\mu \partial^\theta (u_0+w_{m-1})|
+ \sum_{|\mu|\le 90} |Z^\mu [\Box,\eta]u|.
\end{multline}
For the 
first term in the right side of \eqref{zeroCase}, we can apply
Sobolev embedding on $\S^2$ and the Schwarz inequality to obtain
\begin{multline*}
  \sum_{|\mu|\le 90} \int_0^T \|\la x\ra^{-1/2} G_m\|_{L^1_rL^2_\omega}\,ds
\lesssim \varepsilon \\+ 
\sum_{|\nu|\le 48} \|\la x\ra^{-1/4} Z^\nu (u_0+ w_{m-1})\|_{L^2_{t,x}(S_T)}
\sum_{|\mu|\le 90} \|\la x\ra^{-1/4} Z^\mu (u_0+ w_m)'\|_{L^2_{t,x}(S_T)}
\\+ \sum_{|\nu|\le 48} \|\la x\ra^{-1/4} Z^\nu (u_0+w_m)'\|_{L^2_{t,x}(S_T)}
\sum_{\substack{|\mu|\le 90\\|\theta|\le 1}} \|\la x\ra^{-1/4}
Z^\mu \partial^\theta (u_0+w_{m-1})\|_{L^2_{t,x}(S_T)}
\\+ \sum_{|\nu|\le 48} \|\la x\ra^{-1/4} Z^\nu
(u_0+w_{m-1})\|_{L^2_{t,x}(S_T)} \sum_{\substack{|\mu|\le 90\\|\theta|\le
    1}} \|\la x\ra^{-1/4} Z^\mu \partial^\theta (u_0+w_{m-1})\|_{L^2_{t,x}(S_T)}.
\end{multline*}
Here we have also applied \eqref{local}.
The above is, in turn,
\[\la T\ra^{1/2}(M_0(T)+ M_{m-1}(T)) (M_0(T)+M_{m-1}(T)+M_m(T)) +
C_2\varepsilon.\]

The second term in the right side of \eqref{zeroCase} is better
behaved.  Indeed, using $IX$, it follows immediately that this is
\[\lesssim (\log(2+T)) (M_0(T)+M_{m-1}(T))(M_0(T)+M_m(T)).\]

Combining the above, we have established
\[III\Bigr|_{|\nu|=0} + IV\Bigr|_{|\nu|=0}\le
C(M_0(T)+M_{m-1}(T))(M_0(T)+M_{m-1}(T)+M_m(T))\la T\ra^{1/2} +C_2\varepsilon.\]

{\bf \em Bounds for $III$ and $IV$, $|\nu|=1$:}
We, again, fix $\beta$ as above and seek to control the
contribution away from the boundary as the coefficients of $Z$ are
$O(1)$ on the support of $\beta$ and the corresponding terms near the
boundary are controlled by $II$.   Indeed, we apply \eqref{le} and \eqref{DZle} to
$(1-\beta)Z^\mu w_m$, which solves the boundaryless equation shown in
the previous subsection.  This yields
\begin{multline}\label{III_IV_RHS}
  \sup_{t\in [0,T]} \sum_{|\mu|\le 90} \|Z^\mu w_m'(t,\cd)\|_2 +
  \sum_{|\mu|\le 90} \la T\ra^{-1/4} \|\la x\ra^{-1/4} Z^\mu
  w_m'\|_{L^2_{t,x}(S_T)} 
\\\lesssim \sum_{|\mu|\le 90} \int_0^T \|Z^\mu
  \Box w_m(s,\cd)\|_2\,ds + \sum_{|\mu|\le 91} \|\partial^\mu w_m\|_{L^2_{t,x}([0,T]\times \{|x|<3\})}.
\end{multline}
The last term in the above equation is subject to the bounds previously
established for $II$.

Mimicking the arguments above, we obtain to following bound for the
first term in the right side of \eqref{III_IV_RHS}:
\begin{multline*}
  \int_0^T \sum_{|\nu|\le 46} \|Z^\nu
  (u_0+w_{m-1})(s,\cd)\|_\infty \sum_{|\mu|\le 90} \|Z^\mu \partial^2
  (u_0+w_m)(s,\cd)\|_2\,ds
\\+ \int_0^T \sum_{|\nu|\le 46} \|Z^\nu \partial (u_0+w_m)(s,\cd)\|_\infty
\sum_{\substack{|\mu|\le 90\\|\nu|\le 1}} \|Z^\mu \partial^\nu
(u_0+w_{m-1})(s,\cd)\|_2\,ds
\\+ \int_0^T \sum_{|\nu|\le 46} \|Z^\nu (u_0+w_{m-1})(s,\cd)\|_\infty
\sum_{\substack{|\mu|\le 90\\|\nu|\le 1}} \|Z^\mu \partial^\nu
(u_0+w_{m-1})(s,\cd)\|_2\,ds
\\+ \int_0^T \sum_{|\mu|\le 90} \|Z^\mu [\Box,\eta]u(s,\cd)\|_2\,ds.
\end{multline*}
The key thing to note here is that as we are not utilizing estimates
for time-dependent perturbations of $\Box$, we must face a term of the
form $\sum_{|\mu|\le 90} \|Z^\mu w_m''\|_2$, which shall be bounded
in the following section.  Utilizing terms $III$ and $IX$ of \eqref{M} as well as
\eqref{local}, it follows that this piece satisfies
\[III\Bigr|_{|\nu|=1} + IV\Bigr|_{|\nu|=1} \le C(M_0(T)+M_{m-1}(T))(M_0(T)+M_{m-1}(T)+M_m(T))\log(2+T) + C_2\varepsilon.\]

{\em \bf Bounds for $III$, $|\nu|=2$:}
With $\gamma$ chosen as when bounding $I$, we have \eqref{gammaDecay1}
and \eqref{gammaDecay2} with $\delta=C_1 \varepsilon$.
We may, thus, apply \eqref{LZenergy}.  Upon integrating
\eqref{LZenergy} in $t$, applying \eqref{gammaDecay2} and \eqref{lifespan}, and
bootstrapping, we need to control
\begin{equation}\label{III_RHS}\int_0^T \sum_{\substack{|\mu|\le 90\\|\nu|=1}} \|\Box_\gamma
Z^\mu \partial^\nu w_m(t,\cd)\|_2\,dt + \sum_{|\mu|\le 92}
\|\partial^\mu w_m'\|_{L^2_{t,x}([0,T]\times\{|x|<1\})}.\end{equation}
As above, the bound established for $II$ applies to the latter term,
and we need only control the former term.

Using the product rule, it follows that
\begin{multline*}
  \sum_{\substack{|\mu|\le 90\\|\nu|=1}} |\Box_\gamma Z^\mu \partial^\nu w_m|
\lesssim \sum_{|\nu|\le
    47} |Z^\nu (u_0+w_{m-1})| \Bigl[
  \sum_{\substack{|\mu|\le 90\\|\nu|\le 1}} |Z^\mu \partial^\nu (u_0+w_m)'|
+ \sum_{|\mu|\le 93}
|Z^\mu u_0|\Bigr]
\\ 
+ \sum_{|\nu|\le
    47} |Z^\nu (u_0+w_{m})'|
  \sum_{\substack{|\mu|\le 90\\|\nu|\le 2}} |Z^\mu \partial^\nu (u_0+w_{m-1})|
\\+ \sum_{|\nu|\le
    47} |Z^\nu (u_0+w_{m-1})|
  \sum_{\substack{|\mu|\le 90\\|\nu|\le 2}} |Z^\mu \partial^\nu (u_0+w_{m-1})|
+ \sum_{|\mu|\le 91} |Z^\mu [\Box,\eta]u|.
\end{multline*}
And hence, using \eqref{local},
\begin{multline*}
  \int_0^T \sum_{\substack{|\mu|\le 90\\|\nu|=1}} \|\Box_\gamma
Z^\mu \partial^\nu w_m(t,\cd)\|_2\,dt \lesssim\varepsilon \\+ 
\int_0^T \sum_{|\nu|\le
    47} \|Z^\nu (u_0+w_{m-1})(t,\cd)\|_\infty \Bigl[
  \sum_{\substack{|\mu|\le 90\\|\nu|\le 1}} \|Z^\mu \partial^\nu (u_0+w_m)'(t,\cd)\|_2
+ \sum_{|\mu|\le 93}
\|Z^\mu u_0(t,\cd)\|\Bigr]\,dt
\\ 
+ \int_0^T \sum_{|\nu|\le
    47} \|Z^\nu (u_0+w_{m})'(t,\cd)\|_\infty
  \sum_{\substack{|\mu|\le 90\\|\nu|\le 2}} \|Z^\mu \partial^\nu (u_0+w_{m-1})(t,\cd)\|_2\,dt
\\+ \int_0^T \sum_{|\nu|\le
    47} \|Z^\nu (u_0+w_{m-1})(t,\cd)\|_\infty
  \sum_{\substack{|\mu|\le 90\\|\nu|\le 2}} \|Z^\mu \partial^\nu (u_0+w_{m-1})(t,\cd)\|_2\,dt.
\end{multline*}
We now use terms $III$ and $IX$ of \eqref{M}.  This immediately yields
\[III\Bigr|_{|\nu|=2} \le C(M_0(T)+M_{m-1}(T))(M_0(T)+M_{m-1}(T)+M_m(T))\log(2+T) +C_2\varepsilon.\]

{\em \bf Bound for $V$:}
It is here where our approach most differs from that of \cite{DZ}.
When proving energy estimate involving $L$, \cite{DZ} utilized the
star-shapedness, as in \cite{KSS3}, to show that the worst contribution on the boundary
had a favorable sign.  Our approach instead follows that of
\cite{MS3}, which relies instead on \eqref{ms_bdy}.

We first note that by \eqref{ellReg}, it suffices to estimate 
\[\sum_{\substack{j+k\le 81\\k\le 1}} \|L^k \partial_t^j w_m'(t,\cd)\|_2
+ \sum_{\substack{|\mu|+k\le 80\\k\le 1}} \|L^k \partial^\mu \Box w_m(t,\cd)\|_2.\]
We further note that
\[\sum_{\substack{j+k\le 81\\k\le 1}} \|L^k \partial_t^j
w_m'(t,\cd)\|_2 \le \sum_{k\le 80} \|\tilde{L} \partial_t^j
w_m'(t,\cd)\|_2 + \sum_{|\mu|\le 81} \|\partial^\mu w_m'(t,\cd)\|_2.\]
As the latter term is subject to the previously established bounds for
$I$, it suffices to control
\[\sum_{j\le 80} \|\tilde{L} \partial_t^j
w_m'(t,\cd)\|_2  + \sum_{\substack{|\mu|+k\le 80\\k\le 1}}
\|L^k \partial^\mu \Box w_m(t,\cd)\|_2.\]

For the former term, we shall employ \eqref{X} with $\gamma$ as in the
argument for term $I$.  After integrating \eqref{X}, applying
\eqref{gammaDecay2} and \eqref{lifespan}, and bootstrapping, it
remains to bound
\begin{multline}\label{VneedToBound}
 \int_0^T \sum_{\substack{j+k\le 81\\k\le 1}} \Bigl(\|\tilde{L}^k \partial_t^j
    \Box_\gamma w_m(t,\cd)\|_2 
    + \|[\tilde{L}^k\partial_t^j,
    \gamma^{\alpha\beta}\partial_\alpha\partial_\beta]w_m(t,\cd)\|_2\Bigr)\,dt \\+\int_0^T
    \sum_{j\le 80} \|\partial_t^j \Box w_m(t,\cd)\|_2 \,dt +
    \sum_{|\mu|\le 80} \int_0^T \|\partial^\mu
    w_m'(t,\cd)\|_{L^2(|x|<1)}\, dt
\\+
\sup_{t\in [0,T]} \sum_{\substack{|\mu|+k\le 80\\k\le 1}}
\|L^k \partial^\mu \Box w_m(t,\cd)\|_2.
\end{multline}

We note that
\begin{multline}\label{LproductRule}
  \sum_{\substack{j+k\le 81\\k\le 1}} \Bigl(|\tilde{L}^k \partial_t^j
  \Box_\gamma w_m| + |[\tilde{L}^k\partial_t^j, \Box-\Box_\gamma]
  w_m|\Bigr) \lesssim \sum_{|\nu|\le 41} |\partial^\nu
  (u_0+w_{m-1})|
  \sum_{|\mu| \le 80} |L\partial^\mu w_m'|
\\+  \Bigl(\sum_{|\mu|\le 41} (|\partial^\mu (w_{m-1}+u_0)| +
|\partial^\mu (w_m+u_0)'|)\Bigr) \sum_{|\mu|\le 80} |L \partial^\mu
(w_{m-1}+u_0)'|
\\+\sum_{\substack{j+|\mu|\le 41\\j\le 1}} |L^j \partial^\mu
(w_{m-1}+u_0)| \Bigl(\sum_{|\nu|\le 81} (|\partial^\nu (w_m+u_0)'| +
|\partial^\nu(w_{m-1}+u_0)'| + |\partial^\nu u_0|)\Bigr)
\\+\sum_{\substack{j+|\mu|\le 41\\j\le 1}} |L^j \partial^\mu
(w_m+u_0)'| \sum_{|\nu|\le 81} |\partial^\nu (w_{m-1}+u_0)'|
\\+ |w_{m-1}+u_0|^2 + \sum_{|\mu|\le 81} |\partial^\mu [\Box,\eta] u|.
\end{multline}
For the last term, we are using the assumption that the Cauchy data
are compactly supported and finite propagation speed in order to
guarantee that the coefficients of $L$ are $O(1)$ on the support of
$[\Box,\eta]u$.

For the terms on the
third and fourth lines, we shall apply \eqref{wtdSob} on dyadic
intervals and then sum over those dyadic intervals.  See \cite{KSS}
for a more detailed computation.  Upon doing so, we see that the $L^2$
norm of the terms on the third and fourth lines above is bounded by
\begin{multline*}
\sum_{\substack{j+|\mu|\le 43\\j\le 1}} \|\la x\ra^{-1/4} L^j Z^\mu
(w_{m-1}+u_0)\|_2\\\times \Bigl(\sum_{|\nu|\le 81} (\|\la x\ra^{-3/4}\partial^\nu (w_m+u_0)'\|_2 +
\|\la x\ra^{-3/4} \partial^\nu(w_{m-1}+u_0)'\|_2 + \|\la x\ra^{-3/4} \partial^\nu u_0\|_2)\Bigr)
\\+\sum_{\substack{j+|\mu|\le 41\\j\le 1}} \|\la x\ra^{-1/4}L^j Z^\mu
(w_m+u_0)'\|_2 \sum_{|\nu|\le 81} \|\la x\ra^{-3/4}\partial^\nu (w_{m-1}+u_0)'\|_2.
\end{multline*}
When integrated in $t$, we apply the Schwarz inequality and utilize terms
$II$ and $XIII$ of \eqref{M} to establish control.

For the terms in the right side of \eqref{LproductRule} which are on
the first and second lines, we apply $IX$ and $V$ of \eqref{M}.  And
control for the second to last term in \eqref{LproductRule} follows
from terms $IX$ and $III$ of \eqref{M}.

Arguing as such yields the bound
\begin{multline*}
\int_0^T \sum_{\substack{j+k\le 81\\k\le 1}} \Bigl(\|\tilde{L}^k \partial_t^j
    \Box_\gamma w_m(t,\cd)\|_2 
    + \|[\tilde{L}^k\partial_t^j,
    \gamma^{\alpha\beta}\partial_\alpha\partial_\beta]w_m(t,\cd)\|_2\Bigr)\,dt 
\\\le C(M_0(T)+M_{m-1}(T))(M_0(T)+M_{m-1}(T)+M_m(T))\la T\ra^{1/4} + C_2\varepsilon.
  \end{multline*}
The integrand of last term in \eqref{VneedToBound}  is also controlled
by the right side of \eqref{LproductRule}.  As there is no time
integral, it can be bounded much more easily just using Sobolev
embedding and terms $I$ and $V$ of \eqref{M}.  No loss of $\la
T\ra^{1/4}$ is necessitated here.  The third term in
  \eqref{VneedToBound} was previously controlled while establishing
  the bound for, say, $II$.  For this term, as we saw previously,
  logarithmic losses would suffice.

It only remains to establish control for 
\[\sum_{|\mu|\le 80} \int_0^T \|\partial^\mu w_m'(t,\cd)\|_{L^2(|x|<1)}\,dt.\]
We apply \eqref{ms_bdy} and see that it suffices to bound
\[\sum_{|\mu|\le 81} \int_0^T \int_0^s \|\partial^\mu \Box
w_m(\tau,\cd)\|_{L^2(||x|-(s-\tau)|<10)}\,d\tau\,ds.\]
Similar to \eqref{productrule}, we have
\begin{multline*}
\sum_{|\mu|\le 81} |\partial^\mu \Box w_m| \lesssim
\sum_{|\nu|\le
    41} |\partial^\nu (u_0+w_{m-1})|
  \sum_{|\mu|\le 82} |\partial^\mu (u_0+w_m)'|\\ 
+ \sum_{|\nu|\le
    41} |\partial^\nu (u_0+w_{m})'|
  \sum_{|\mu|\le 81} |\partial^\mu (u_0+w_{m-1})'|
+ \sum_{|\nu|\le
    41} |\partial^\nu (u_0+w_{m-1})|
  \sum_{|\mu|\le 81} |\partial^\mu (u_0+w_{m-1})'|
\\+ |u_0+w_{m-1}|^2
+ \sum_{|\mu|\le 81} |\partial^\mu [\Box,\eta]u|.
\end{multline*}
With the exception of the last term above, for which we instead use
\eqref{local}, we apply \eqref{wtdSob2} to see that
\begin{multline*}
\sum_{|\mu|\le 81} \|\partial^\mu \Box w_m\|_{L^2(||x|-(s-\tau)|<10)} \\\lesssim
\sum_{|\nu|\le
    43} \|\la x\ra^{-1/4} Z^\nu (u_0+w_{m-1})\|_{L^2(||x|-(s-\tau)|<20)}
  \sum_{|\mu|\le 82} \|\la x\ra^{-3/4} \partial^\mu (u_0+w_m)'\|_{L^2(||x|-(s-\tau)|<20)}\\ 
+ \sum_{|\nu|\le
    43} \|\la x\ra^{-1/4}Z^\nu (u_0+w_{m})'\|_{L^2(||x|-(s-\tau)|<20)}
  \sum_{|\mu|\le 81} \|\la x\ra^{-3/4} \partial^\mu (u_0+w_{m-1})'\|_{L^2(||x|-(s-\tau)|<20)}
\\+ \sum_{|\nu|\le
    43} \|\la x\ra^{-1/4} Z^\nu (u_0+w_{m-1})\|_{L^2(||x|-(s-\tau)|<20)}
  \sum_{|\mu|\le 81} \|\la x\ra^{-3/4}\partial^\mu (u_0+w_{m-1})'\|_{L^2(||x|-(s-\tau)|<20)}
\\+ \sum_{|\nu|\le 2} \|\la x\ra^{-1/4}  Z^\nu(u_0+w_{m-1})\|^2_{L^2(||x|-(s-\tau)|<20)}
+ \sum_{|\mu|\le 81} \|\partial^\mu [\Box,\eta]u\|_{L^2(||x|-(s-\tau)|<10)}.
\end{multline*}
Since the sets $\{||x|-(j-\tau)|<20\}$ have finite overlap as $j$
ranges over the nonnegative integers, it follows that upon integrating
in $\tau$ and $s$ that
\begin{multline*}
\sum_{|\mu|\le 81} \int_0^T \int_0^s \|\partial^\mu \Box w_m\|_{L^2(||x|-(s-\tau)|<10)} \,d\tau\,ds\\\lesssim
\sum_{|\nu|\le
    43} \|\la x\ra^{-1/4} Z^\nu (u_0+w_{m-1})\|_{L^2_{t,x}(S_T)}
  \sum_{|\mu|\le 82} \|\la x\ra^{-3/4} \partial^\mu (u_0+w_m)'\|_{L^2_{t,x}(S_T)}\\ 
+ \sum_{|\nu|\le
    43} \|\la x\ra^{-1/4}Z^\nu (u_0+w_{m})'\|_{L^2_{t,x}(S_T)}
  \sum_{|\mu|\le 81} \|\la x\ra^{-3/4} \partial^\mu (u_0+w_{m-1})'\|_{L^2_{t,x}(S_T)}
\\+ \sum_{|\nu|\le
    43} \|\la x\ra^{-1/4} Z^\nu (u_0+w_{m-1})\|_{L^2_{t,x}(S_T)}
  \sum_{|\mu|\le 81} \|\la x\ra^{-3/4}\partial^\mu (u_0+w_{m-1})'\|_{L^2_{t,x}(S_T)}
\\+ \sum_{|\nu|\le 2} \|\la x\ra^{-1/4}  Z^\nu(u_0+w_{m-1})\|^2_{L^2_{t,x}(S_T)}
+\varepsilon.
\end{multline*}
Here we have also employed \eqref{local}.  And thus, we see that this
boundary term is 
\[\le C(M_0(T)+M_{m-1}(T))(M_0(T)+M_{m-1}(T)+M_{m}(T))\la T\ra^{1/2} + C_2\varepsilon,\]
and hence, when combined with the above
\[V\le C(M_0(T)+M_{m-1}(T))(M_0(T)+M_{m-1}(T)+M_{m}(T))\la T\ra^{1/2} + C_2\varepsilon.\]

{\em \bf Bound for $VI$:}  The arguments to bound terms $VI$, $VII$,
and $VIII$ follow those of the corresponding terms with no $L$ quite
closely.  Indeed, for $VI$, we, as in the bound for $II$, apply
\eqref{DZle} to $L \partial^\mu w_m$ cutoff away from the boundary and
\eqref{locEnergyExt} to both the remaining portion as well as the
compactly supported commutator which results from cutting off above.
It remains then to control
\begin{equation}\label{VI_RHS}\sum_{\substack{j+|\mu|\le 77\\j\le 1}} \int_0^T \|L^j\partial^\mu
\Box w_m(s,\cd)\|_2\,ds + \sum_{\substack{j+|\mu|\le 75\\j\le 1}}
\|L^j \partial^\mu \Box w_m\|_{L^2_{t,x}(S_T)}.\end{equation}
By Sobolev embeddings, it suffices to control the first term.

We must now take a little care with the location of the scaling vector
field.  Playing the role of \eqref{productrule}, we have
\begin{multline*}
\sum_{\substack{j+|\mu|\le 77\\j\le 1}} |L^j \partial^\mu \Box w_m| \lesssim
\sum_{|\nu|\le
    40} |\partial^\nu (u_0+w_{m-1})|
  \sum_{\substack{j+|\mu|\le 78\\j\le 1}} |L^j\partial^\mu (u_0+w_m)'|\\
\sum_{\substack{j+|\nu|\le
    40\\j\le 1}} |L^j \partial^\nu (u_0+w_{m-1})|
  \sum_{|\mu|\le 78} |\partial^\mu (u_0+w_m)'|\\ 
+ \sum_{|\nu|\le
    41} |\partial^\nu (u_0+w_{m})'|
  \sum_{\substack{j+|\mu|\le 77\\j\le 1}} |L^j \partial^\mu
  (u_0+w_{m-1})'|
\\+ \sum_{\substack{j+|\nu|\le
    41\\j\le 1}} |L^j \partial^\nu (u_0+w_{m})'|
  \sum_{|\mu|\le 77} |\partial^\mu (u_0+w_{m-1})'|
\\+ \sum_{|\nu|\le
    40} |\partial^\nu (u_0+w_{m-1})|
  \sum_{\substack{j+|\mu|\le 77\\j\le 1}} |L^j \partial^\mu (u_0+w_{m-1})'|
\\+ \sum_{\substack{j+|\nu|\le
    40\\j\le 1}} |L^j \partial^\nu (u_0+w_{m-1})|
  \sum_{|\mu|\le 77} |\partial^\mu (u_0+w_{m-1})'|
\\+ |u_0+w_{m-1}| \sum_{j\le 1} |L^j (u_0+w_{m-1})|
+ \sum_{\substack{j+|\mu|\le 78\\j\le 1}} |L^j \partial^\mu [\Box,\eta]u|.
\end{multline*}
When the scaling vector field is on the higher order term, we shall
generally utilize terms $V$ and $IX$ of \eqref{M}.  Alternatively,
when the scaling vector field is on the lower order term, we utilize
\eqref{wtdSob}.  More specifically, we decompose dyadically in $x$,
apply \eqref{wtdSob}, and apply the Schwarz inequality in both the
dyadic summation variable and in $t$ as above.  Upon doing so, we
can bound utilizing $II$, $IV$ and $VII$ instead.  For the second to last
term, $VII$ and $IX$ are employed.  And finally, \eqref{local}
provides the bound for the final term. 

We illustrate arguing in this fashion by examining the
$L^1([0,T];L^2(\ext))$-norm of the first two terms in the right.  Indeed,
the norm of these terms is bounded by
\begin{multline*}
  \int_0^T \sum_{|\nu|\le
    40} \|\partial^\nu (u_0+w_{m-1})(s,\cd)\|_\infty
  \sum_{\substack{j+|\mu|\le 78\\j\le 1}} \|L^j\partial^\mu (u_0+w_m)'(s\cd)\|_2\,ds\\
\sum_{\substack{j+|\nu|\le
    42\\j\le 1}} \|\la x\ra^{-1/4} L^j Z^\nu (u_0+w_{m-1})\|_{L^2_{t,x}(S_T)}
  \sum_{|\mu|\le 78} \|\la x\ra^{-3/4} \partial^\mu (u_0+w_m)'\|_{L^2_{t,x}(S_T)}.
\end{multline*}
And using $XI$, $V$, $VIII$, and $II$, this is
\[\lesssim (M_0(T)+M_{m-1}(T))(M_0(T)+M_{m}(T))\la T\ra^{1/4}.\]
The remaining terms above are handled in a directly analogous way,
which yields
\[
VI\le C
(M_0(T)+M_{m-1}(T))(M_0(T)+M_{m-1}(T)+M_m(T))\la T\ra^{1/4} + C_2\varepsilon.
\]

{\em \bf Bound for $VII$ and $VIII$ with $|\nu|=0$:}
We fix $\beta$ as above.  It suffices to bound $(1-\beta)L Z^\mu w_m$
as the Dirichlet boundary conditions and the boundedness of the
coefficients of $Z$ on the support of $\beta$ allow us to control
$\beta L Z^\mu w_m$ using term $VI$.
For $\mu$ fixed, we decompose
\begin{multline*}
(1-\beta)L Z^\mu \Box w_m
\\= \partial_\alpha[(1-\eta)(1-\beta) (b^{\alpha\beta}
(u_0+w_{m-1})  L Z^\mu \partial_\beta (u_0+w_m) +
b^{\alpha\beta}_\gamma \partial_\gamma (u_0+w_{m-1})  L
Z^\mu \partial_\beta (u_0+w_m))] 
\\+\tilde{G}_m(t,x).
\end{multline*}
We apply \eqref{divForm}, \eqref{noD}, and \eqref{noDcomm}.  The
latter is used for $[\Box,\beta] LZ^\mu w_m$.  This yields
\begin{multline}\label{LzeroCase}
  \sup_{t\in [0,T]} \sum_{|\mu|\le 70} \|(1-\beta) L Z^\mu
  w_m(t,\cd)\|_2 + \sum_{|\mu|\le 70} \la T\ra^{-1/4} \|\la
  x\ra^{-1/4} (1-\beta)L Z^\mu w_m\|_{L^2_{t,x}(S_T)}
\\\lesssim \int_0^T
\|\la x\ra^{-1/2}\tilde{G}_m(s,\cd)\|_{L^1_rL^2_\omega}\,ds 
+\sum_{|\mu|\le 70}\sum_{|\theta|\le 1}\int_0^T \|
|\partial^\theta (u_0+w_{m-1})| |L Z^\mu (u_0+w_m)'|\|_2\,ds
\\+ \sum_{|\mu|\le 70} \|L \partial^\mu w_m'\|_{L^2_{t,x}([0,T]\times\{|x|<3\})}.
\end{multline}
The last term is controlled by $VI$, and the bound previously
established for $VI$ shall simply be cited to control this piece.

To control the first term in the right of \eqref{LzeroCase}, we need
not be as precise with the location of the scaling vector fields, and
indeed, we can crudely utilize the obvious analog of
\eqref{productrule2} where every term on the right is permitted at
most one occurrence of $L$.  Upon doing so and using Sobolev embedding
on $\S^2$ as well as \eqref{local}, we obtain
\begin{multline*}
  \sum_{|\mu|\le 70} \int_0^T \|\la x\ra^{-1/2} \tilde{G}_m\|_{L^1_rL^2_\omega}\,ds
\lesssim \varepsilon \\+ 
\sum_{\substack{|\nu|\le 38\\ j\le 1}} \|\la x\ra^{-1/4} L^j Z^\nu (u_0+w_{m-1})\|_{L^2_{t,x}(S_T)}
\sum_{\substack{|\mu|\le 70\\k\le 1}} \|\la x\ra^{-1/4} L^k Z^\mu (u_0+w_m)'\|_{L^2_{t,x}(S_T)}
\\+ \sum_{\substack{|\nu|\le 38\\j\le 1}} \|\la x\ra^{-1/4} L^j Z^\nu (u_0+w_m)'\|_{L^2_{t,x}(S_T)}
\sum_{\substack{|\mu|\le 70\\k\le 1\\|\theta|\le 1}} \|\la x\ra^{-1/4}
L^k Z^\mu \partial^\theta (u_0+w_{m-1})\|_{L^2_{t,x}(S_T)}
\\+ \sum_{\substack{|\nu|\le 38\\j\le 1}} \|\la x\ra^{-1/4} L^j Z^\nu
(u_0+w_{m-1})\|_{L^2_{t,x}(S_T)} \sum_{\substack{|\mu|\le 70\\k\le 1\\|\theta|\le
    1}} \|\la x\ra^{-1/4}L^k Z^\mu \partial^\theta (u_0+w_{m-1})\|_{L^2_{t,x}(S_T)}.
\end{multline*}
Each of these individual factors is controlled either by $\la
T\ra^{1/4} IV$ or $\la T\ra^{1/4} VIII$.  Thus, the first term in the
right of \eqref{LzeroCase} is
\[\le C\la T\ra^{1/2} (M_0(T)+M_{m-1}(T))(M_0(T)+M_{m-1}(T)+M_m(T)) + C_2\varepsilon.\]

For the second term in the right side of \eqref{LzeroCase}, we may
simply use $IX$ and $VII$ to immediately see that it is
\[\le C (M_0(T)+M_{m-1}(T))(M_0(T)+M_m(T))\log(2+T).\]
And thus, by combining the above bounds, we see that
\[VII\Bigr|_{|\nu|=0} + VIII\Bigr|_{|\nu|=0} \le C\la T\ra^{1/2}
(M_0(T)+M_{m-1}(T))(M_0(T)+M_{m-1}(T)+M_m(T)) + C_2\varepsilon.\]

{\em  \bf Bound for $VII$ and $VIII$ with $|\nu|=1$:}
Here, again, it suffices to control the given norm when $w_m$ is
replaced by $(1-\beta) w_m$.  As the coefficients of $Z$ are $O(1)$ on
the support of $\beta$, the corresponding $\beta w_m$ terms are
subject to the bounds previously established for $V$ and $VI$.

Applying \eqref{le} and \eqref{DZle} to $(1-\beta) L Z^\mu w_m$, one
obtains 
\begin{multline*}
  \sup_{t\in [0,T]} \sum_{|\mu|\le 70} \|L Z^\mu w_m'(t,\cd)\|_2 +
  \sum_{|\mu|\le 70} \la T\ra^{-1/4} \|\la x\ra^{-1/4} L Z^\mu
  w_m'\|_{L^2_{t,x}(S_T)} 
\\\lesssim \sum_{\substack{|\mu|\le 70\\j\le 1}} \int_0^T \|L^j Z^\mu
  \Box w_m(s,\cd)\|_2\,ds + \sum_{|\mu|\le 71} \|L \partial^\mu w_m\|_{L^2_{t,x}([0,T]\times \{|x|<3\})}.
\end{multline*}
The bound for term $VI$ applies to the latter term on the right side.

Here we need to pay attention to the location of $L$ for terms
involving $\partial^2 w_m$, but for the other terms we may be more crude and simply
permit up to one occurrence of $L$ on each factor.  For everything except
for the case that all of the vector fields land on $\partial^2 w_m$,
we dyadically decompose, apply \eqref{wtdSob}, and use the
Cauchy-Schwarz inequality in the dyadic variable as well as $t$.  
Upon
doing so, we obtain that first term in the right side of the preceding
equation is controlled by
\begin{multline*}
  \int_0^T \sum_{|\nu|\le 36} \|Z^\nu
  (u_0+w_{m-1})(s,\cd)\|_\infty \sum_{|\mu|\le 70} \|L Z^\mu \partial^2
  (u_0+w_m)(s,\cd)\|_2\,ds
\\
+ \sum_{\substack{|\nu|\le 38\\j\le 1}} \|\la x\ra^{-1/2} L^j Z^\nu
(u_0+w_{m-1})\|_{L^2_{t,x}(S_T)} \sum_{|\mu|\le 71} \|\la x\ra^{-1/2}
Z^\mu \partial (u_0+w_m)\|_{L^2_{t,x}(S_T)}
\\+ \sum_{\substack{|\nu|\le 38\\j\le 1}} \|\la x\ra^{-1/2} L^j Z^\nu \partial (u_0+w_m)\|_{L^2_{t,x}(S_T)}
\sum_{\substack{|\mu|\le 70\\k\le 1\\|\nu|\le 1}} \|\la x\ra^{-1/2}
L^k Z^\mu \partial^\nu
(u_0+w_{m-1})\|_{L^2_{t,x}(S_T)}
\\+ \sum_{\substack{|\nu|\le 38\\j\le 1}} \|\la x\ra^{-1/2} L^j Z^\nu (u_0+w_{m-1})\|_{L^2_{t,x}(S_T)}
\sum_{\substack{|\mu|\le 70\\k\le 1\\|\nu|\le 1}} \|\la x\ra^{-1/2}
L^k Z^\mu \partial^\nu
(u_0+w_{m-1})\|_{L^2_{t,x}(S_T)}
\\+ \int_0^T \sum_{\substack{|\mu|\le 70\\j\le 1}} \|L^j Z^\mu [\Box,\eta]u(s,\cd)\|_2\,ds.
\end{multline*}
The first term is bounded using $VII$ and $IX$, while for the
next three, term $VIII$ is primarily used.  The only exception is the
second factor in the second term where $IV$ is applied.  The estimate
\eqref{local} is used to control the last term.  This shows that the
above yields
\[VII\Bigr|_{|\nu|=1} + VIII\Bigr|_{|\nu|=1} \le C(M_0(T)+M_{m-1}(T))(M_0(T)+M_{m-1}(T)+M_m(T))\la T\ra^{1/2} +
C_2\varepsilon.\] 
Because the weights that appear in the $L^2_{t,x}$ norms are $\la
x\ra^{-1/2}$ rather than $\la x\ra^{-1/4}$, it would be rather easy to
replace $\la T\ra^{1/2}$ in this bound by $\log(2+T)$.  As this will
not improve the final lifespan, we choose to not lengthen the argument
in order to show this.

{\em\bf Bound for $VII$ with $|\nu|=2$:}
With $\gamma$ chosen as in \eqref{gamma}, which satisfies
\eqref{gammaDecay1} and \eqref{gammaDecay2} with
$\delta=C_1\varepsilon$, we employ \eqref{LZenergy}.  After
integrating, applying \eqref{gammaDecay2}, using that $T\le
T_\varepsilon$ and bootstrapping, it suffices to bound
\[\int_0^T \sum_{\substack{|\mu|\le 70\\|\nu|=1}} \|\Box_\gamma L
Z^\mu \partial^\nu w_m(t,\cd)\|_2\,dt + \sum_{\substack{|\mu|\le
    72\\k\le 1}} \|L^k \partial^\mu w_m'\|_{L^2_{t,x}([0,T]\times
  \{|x|<1\})}.\]
The latter term is controlled by $VI$ of \eqref{M}, and the bound
proved previously for that term suffices.

Here, we have
\begin{multline*}
  \sum_{\substack{|\mu|\le 70\\|\nu|=1}} |\Box_\gamma L
  Z^\mu \partial^\nu w_m|\lesssim 
\sum_{|\nu|\le 37} |Z^\nu (u_0+w_{m-1})|
\Bigl[\sum_{\substack{|\mu|\le 70\\|\nu|\le 1\\j\le 1}} |L^j
Z^\mu \partial^\nu (u_0+w_m)'| + \sum_{|\mu|\le 73} |\partial^\mu u_0|\Bigr]
\\+ \sum_{\substack{|\nu|\le 37\\j\le 1}} |L^j Z^\nu (u_0+w_{m-1})|
\Bigl[\sum_{|\mu|\le 71} |Z^\mu (u_0+w_m)'| + \sum_{|\mu|\le 73} |\partial^\mu u_0|\Bigr]
\\+ \sum_{|\nu|\le 37} |Z^\nu (u_0+w_m)'| \sum_{\substack{|\mu\le
    70\\|\nu|\le 2\\j\le 1}} |L^j Z^\mu \partial^\nu (u_0+w_{m-1})|
\\+ \sum_{\substack{|\nu|\le 37\\j\le 1}} |L^j Z^\nu (u_0+w_m)'| \sum_{|\mu\le
    72} |Z^\mu (u_0+w_{m-1})|
\\+ \sum_{|\nu|\le 37} |Z^\nu (u_0+w_{m-1})| \sum_{\substack{|\mu\le
    70\\|\nu|\le 2\\j\le 1}} |L^j Z^\mu \partial^\nu (u_0+w_{m-1})|
\\+ \sum_{\substack{|\nu|\le 37\\j\le 1}} |L^j Z^\nu (u_0+w_{m-1})| \sum_{|\mu\le
    72} |Z^\mu (u_0+w_{m-1})|
+\sum_{|\mu|\le 72} |\partial^\mu [\Box,\eta]u|.
\end{multline*}
Here we have used the assumption that the Cauchy data are compactly
supported and finite propagation speed guarantees that the coefficients
of $L$ and $Z$ are $O(1)$ on the supports of $u_0$ and $[\Box,\eta]u$.

The method of bounding the above terms in $L^1_t([0,T]; L^2_x)$
depends on the location of the scaling vector field.  We will
illustrate the method on the terms in the third and fourth lines
above.  The remaining pieces are controlled in a directly analogous
manner.  When the scaling vector field is on the higher order factor,
we shall use $IX$ and $VII$ (and $III$ in the case that no $L$
appears).  When $L$ lands on the lower order piece, we shall instead
apply \eqref{wtdSob} as above.  Doing so gives the following upper
bound on the $L^1_t([0,T];L^2_x)$ norm of the terms in the third and
fourth lines above:
\begin{multline*}
   \int_0^T \sum_{|\nu|\le 37} \|Z^\nu (u_0+w_m)'(t,\cd)\|_\infty \sum_{\substack{|\mu\le
    70\\|\nu|\le 2\\j\le 1}} \|L^j Z^\mu \partial^\nu (u_0+w_{m-1})(t,\cd)\|_2\,dt
\\+ \sum_{\substack{|\nu|\le 39\\j\le 1}} \|\la x\ra^{-1/2} L^j Z^\nu (u_0+w_m)'\|_{L^2_{t,x}(S_T)} \sum_{|\mu\le
    72} \|\la x\ra^{-1/2} Z^\mu (u_0+w_{m-1})\|_{L^2_{t,x}(S_T)}.
\end{multline*}
By now citing $IX$, $VII$, and $VIII$ of \eqref{M} and arguing
similarly for the remaining nonlinear terms, we see that
\[VII\Bigr|_{|\nu|=2} \le
C(M_0(T)+M_{m-1}(T))(M_0(T)+M_{m-1}(T)+M_m(T)) \la T\ra^{1/2} + C_2\varepsilon.\]

{\em \bf Bound for $IX$:}
Here, we apply \eqref{l1linf}, and we are left with bounding
\begin{equation}\label{RHSl1linf}\int_0^T \int \sum_{\substack{|\mu|+j\le 67\\j\le 1}} |L^j Z^\mu
\Box w_m(s,y)|\,\frac{dy\,ds}{|y|} + \int_0^T
\sum_{\substack{|\mu|+j\le 64\\j\le 1}} \|L^j \partial^\mu \Box
w_m(s,\cd)\|_{L^2(|x|<2)}\,ds.\end{equation}
For the first term, we use an analog of \eqref{productrule} for the
given vector fields and apply the Schwarz inequality.  Upon doing so,
we have
\begin{multline*}
 \sum_{\substack{|\mu|+j\le 69\\j\le 1}} \|\la x\ra^{-1/2} L^j Z^\mu
 (u_0+w_m)\|_{L^2_{t,x}(S_T)} \sum_{\substack{|\mu|+j\le 68\\j\le 1}} \|\la
 x\ra^{-1/2} L^j Z^\mu (u_0+w_{m-1})\|_{L^2_{t,x}(S_T)}
\\+
\Bigl(\sum_{\substack{|\mu|+j\le 68\\j\le 1}} \|\la x\ra^{-1/2} L^j
Z^\mu (u_0+w_{m-1})\|_{L^2_{t,x}(S_T)}\Bigr)^2
+ C_2\varepsilon.
\end{multline*}
Control for this term follows from $IV$ and $VIII$ of
\eqref{M}.  And the second term of \eqref{RHSl1linf} was previously
controlled in the process of bounding term $VI$ of \eqref{M}.  It
follows that
\[IX \le C(M_0(T)+M_{m-1}(T))(M_0(T)+M_{m-1}(T)+M_m(T))\la T\ra^{1/2}
+ C_2\varepsilon.\]

{\em \bf Boundedness of $M_m(T)$:}
If we combine the estimates for terms $I$, \dots, $IX$ just
established, it follows that
\[M_m(T)\le C(M_0(T)+M_{m-1}(T))(M_0(T)+M_{m-1}(T)+M_m(T))\la
T\ra^{1/2} + C_2\varepsilon\]
where $C_2$ is now a fixed constant which is independent of $m$,
$\varepsilon$, and $T$.  If $C_1$ is chosen so that $C_1>2C_2$ and if
we apply \eqref{basecase} as well as the inductive hypothesis, it
follows that
\[M_m(T)\le C\varepsilon (\varepsilon +M_m(T))\la T\ra^{1/2} + \frac{C_1}{2}\varepsilon.\]
If we use \eqref{lifespan}, then if $c$ in \eqref{lifespan} is
sufficiently small we may bootstrap in such a way that
we obtain \eqref{bddness}
if $\varepsilon$ is small enough,
as desired.

{\em \bf Convergence of the sequence $\{w_m\}$:} We now complete the
proof of Theorem \ref{thm1} by showing that the sequence $\{w_m\}$
is Cauchy.  Standard results show that the limiting function solves
\eqref{reduced}, which is equivalent to the existence promised in
Theorem \ref{thm1}.

Indeed, if we set
\[A_m(T) = \sup_{t\in [0,T]} \sum_{|\mu|\le 90} \|\partial^\mu (w_m-w_{m-1})(t,\cd)\|_2,\]
we may argue quite similarly to the above.  Upon doing so and using
\eqref{bddness}, it can be shown that
\[A_m(T)\le \frac{1}{2}A_{m-1}(T)\]
for $T\le T_\varepsilon$, which completes the proof.



\begin{thebibliography}{MA}



 


\bibitem{BZ} N. Burq and M. Zworski: {\em Geometric control in the
  presence of a black box}.  J. Amer. Math. Soc. {\bf 17} (2004), 443--471.


\bibitem{C} H. Christianson:  {\em Dispersive estimates for manifolds
  with one trapped orbit}.  Comm. Partial Differential Equations {\bf
  33} (2008), 1147--1174.

\bibitem{CW} H. Christianson and J. Wunsch: {\em Local smoothing for
  the Schr\"odinger equation with a prescribed loss}, preprint.
  (ArXiv: 1103.3908)

\bibitem{CdVP} Y. Colin de Verdi\'ere and B. Parisse: {\em \'Equilibre
  instable en r\'egime semi-classique. I. Concentration microlocale}.
  Comm. Partial Differential Equations {\bf 19} (1994), 1535--1563.

\bibitem{DMSZ}Y. Du, J. Metcalfe, C. D. Sogge and Y. Zhou: 
{\em Concerning the Strauss conjecture and almost global existence 
for nonlinear Dirichlet-wave equations in $4$-dimensions}.
Comm. Partial Differential Equations {\bf 33} (2008), 1487--1506.

\bibitem{DZ} Y. Du and Y. Zhou: 
{\em The life span for nonlinear wave equation outside of star-shaped
obstacle in three space dimensions}. Comm. Partial Differential
Equations {\bf 33} (2008), 1455--1486.








\bibitem{HY} K. Hidano and K. Yokoyama: {\em A remark on the almost
  global existence theorems of Keel, Smith, and Sogge}.
  Funkcial. Ekvac.  {\bf 48} (2005), 1--34.

\bibitem{Hormander} L. H\"ormander: {\em On the fully nonlinear Cauchy
  problem with small data. II.}  Microlocal analysis and nonlinear
  waves (Minneapolis, MN, 1988--1989), 51--81, IMA Vol. Math. Appl.,
  30, Springer, New York, 1991.

\bibitem{Htext} L. H\"ormander.  Lectures on nonlinear hyperbolic
  differential equations.  Springer-Verlag, Berlin, 1997.

\bibitem{HormanderL1Linfty} L. H\"ormander: {\em $L^1,\ L^\infty$
  estimates for the wave operator}.  Analyse math\'ematique et
  applications, 211--234, Gauthier-Villars, Montrouge, 1988.


\bibitem{Ikawa1} M. Ikawa: {\em Decay of solutions of the wave
  equation in the exterior of two convex bodies}.  Osaka J. Math. {\bf
  19} (1982), 459--509.

\bibitem{Ikawa2} M. Ikawa: {\em Decay of solutions of the wave
  equation in the exterior of several convex bodies}.
  Ann. Inst. Fourier (Grenoble) {\bf 38} (1988), 113--146.


\bibitem{JK} F. John and S. Klainerman: {\em Almost global existence
  to nonlinear wave equations in three space dimensions}.  Comm. Pure
  Appl. Math. {\bf 37} (1984), 443--455.

\bibitem{KK} S. Katayama and H. Kubo: {\em An alternative proof of
  global existence for nonlinear wave equations in an exterior
  domain}.  J. Math. Soc. Japan {\bf 60} (2008), 1135--1170.

\bibitem{KSS4} M. Keel, H. F. Smith, and C. D. Sogge: {\em On global
  existence for nonlinear wave equations outside of convex
  obstacles}.  Amer. J. Math. {\bf 122} (2000), 805--842.

\bibitem{KSS}  M. Keel, H. F. Smith, and C. D. Sogge: 
{\em Almost global existence for some semilinear wave equations}, 
J. Anal. Math. {\bf 87} (2002), 265-279.

\bibitem{KSS2} M. Keel, H. F. Smith, and C. D. Sogge: {\em Global
  existence for a quasilinear wave equation outside of star-shaped
  domains}.  J. Funct. Anal. {\bf 189} (2002), 155--226.

\bibitem{KSS3} M. Keel, H. F. Smith, and C. D. Sogge: {\em Almost
  global existence for quasilinear wave equations in three space
  dimensions}.  J. Amer. Math. Soc. {\bf 17} (2004), 109--153.


\bibitem{KPV} C. E. Kenig, G. Ponce, and L. Vega: {\em On the Zakharov
  and Zakharov-Schulman systems}.  J. Funct. Anal. {\bf 127} (1995), 204--234.

\bibitem{Klainerman2} S. Klainerman: {\em Uniform decay estimates and
  the Lorentz invariance of the classical wave equation}.  Comm. Pure
  Appl. Math. {\bf 38} (1985), 321--332.

\bibitem{Klainerman} S. Klainerman: {\em The null condition and global
  existence to nonlinear wave equations}.  Lect. Appl. Math. {\bf 23}
  (1986), 293--326.



\bibitem{LMP} P. D. Lax, C. S. Morawetz, and R. S. Phillips: {\em
  Exponential decay of solutions of the wave equation in the exterior
  of a star-shaped obstacle}.  Comm. Pure Appl. Math. {\bf 16} (1963), 477-486.



\bibitem{Lindblad} H. Lindblad: {\em On the lifespan of solutions of
  nonlinear wave equations with small initial data}.  Comm. Pure
  Appl. Math. {\bf 43} (1990), 445--472.






\bibitem{MNS} J. Metcalfe, M. Nakamura, and C. D. Sogge: {\em Global
  existence of solutions to multiple speed systems of quasilinear wave
  equations in exterior domains}.  Forum Math. {\bf 17} (2005), 133--168.

\bibitem{MNS2} J. Metcalfe, M. Nakamura, and C. D. Sogge: {\em Global
  existence of quasilinear, nonrelativistic wave equations satisfying
  the null condition}.  Japan. J. Math. {\bf 21} (2005), 391--472.

\bibitem{MS3} J. Metcalfe and C. D. Sogge: {\em Hyperbolic trapped
  rays and global existence of quasilinear wave equations}.
  Invent. Math. {\bf 159} (2005), 75--117.

\bibitem{MS} J. Metcalfe and C. D. Sogge: {\em Long time existence of
  quasilinear wave equations exterior to star-shaped obstacles via
  energy methods}.  SIAM J. Math. Anal. {\bf 38} (2006), 188--209.


\bibitem{MS4} J. Metcalfe and C. D. Sogge: {\em Global existence of
  null-form wave equations in exterior domains}.  Math. Z. {\bf 256}
  (2007), 521--549.

\bibitem{MS5} J. Metcalfe and C. D. Sogge: {\em Global existence for
    high dimensional quasilinear wave equations exterior to
    star-shaped obstacles}.  Discrete Cts. Dyn. Sys. {\bf 28} (2010), 1589--1601.



\bibitem{MT1} J. Metcalfe and D. Tataru: {\em Global parametrices and
    dispersive estimates for variable coefficient wave equations}.
  Math. Ann., to appear.  (arXiv: 0707.1191)

\bibitem{MT2} J. Metcalfe and D. Tataru: {\em Decay estimates for
    variable coefficient wave equations in exterior domains}.
  Advances in Phase Space Analysis of Partial Differential Equations,
  In Honor of Ferruccio Colombini's 60th Birthday.  Progress in
  Nonlinear Differential Equations and Their Applications, Vol. 78,
  2009, p. 201--217.

\bibitem{Morawetz} C. S. Morawetz: 
{\em The decay of solutions of the exterior initial-boundary value
  problem for the wave equation}.  Comm. Pure Appl. Math. {\bf 14}
(1961), 561--568.

\bibitem{Morawetz2} C. S. Morawetz: {\em Time decay for the nonlinear
  Klein-Gordon equations}.  Proc. Roy. Soc. Ser. A {\bf 306} (1968), 291--296.

\bibitem{MRS} C. S. Morawetz, J. Ralston and W. Strauss: 
{\em Decay of solutions of the wave equation outside nontrapping obstacles}.
Comm. Pure Appl. Math. {\bf 30} (1977), 87--133.


\bibitem{NZ} S. Nonnenmacher and M. Zworski: {\em Semiclassical
  resolvent estimates in chaotic scattering}.
  Appl. Math. Res. Express. {\bf 1} (2009), 74--86.



\bibitem{Ralston} J. V. Ralston: {\em Solutions of the wave equation
  with localized energy}.  Comm. Pure Appl. Math. {\bf 22} (1969), 807--823.




\bibitem{SS} H. F. Smith and C. D. Sogge: 
{\em Global Strichartz estimates for nontrapping perturbations of 
the Laplacian}. 
Comm. Partial Differential Equations {\bf 25} (2000), 2171--2183.



\bibitem{So2} C. D. Sogge.
Lectures on nonlinear wave equations, 2nd edition, 
International Press, Boston, MA, 2008.

\bibitem{Sterb} J. Sterbenz: {\em Angular regularity and Strichartz
  estimates for the wave equation.  With an appendix by Igor
  Rodnianski}.  Int. Math. Res. Not. {\bf 2005}, 187--231.

\bibitem{Strauss} W. A. Strauss: {\em Dispersal of waves vanishing on
  the boundary of an exterior domain}.  Comm. Pure Appl. Math. {\bf
  28} (1975), 265--278.




\bibitem{WZ} J. Wunsch and M. Zworski: {\em Resolvent estimates for
  normally hyperbolic trapped sets}, preprint.  (ArXiv: 1003.4640)



\end{thebibliography}
\end{document}